\documentclass{amsart} %putting equation numbers on the left side

\usepackage{amsfonts,amsmath,amsthm, graphics}

\newtheorem{theorem}{Theorem}%[section]
\newtheorem{lemma}{Lemma}%[section]
\DeclareRobustCommand*\cal{\relax\mathcal}
\newtheorem{definition}{Definition}%[section]
\newtheorem{corollary}{Corollary}%[section]
\newtheorem{remark}{Remark}%[section]
\def\ben{\begin{eqnarray}}
\def\eeqn{\end{eqnarray}}

\def\bd S{\cal C}

\def\>{\geq}

\def\bd{\partial}
\def\bege{\begin{equation}} \def\ende{\end{equation}}

\def\begr{\begin{eqnarray}} \def\endr{\end{eqnarray}}

\def\bege{\begin{equation}} \def\ende{\end{equation}}
\def\begr{\begin{eqnarray}} \def\endr{\end{eqnarray}}
\def\bnum{\begin{enumerate}} \def\enum{\end{enumerate}}
\def\bege{\begin{equation}} \def\ende{\end{equation}}
\def\begr{\begin{eqnarray}} \def\endr{\end{eqnarray}}

\usepackage{amsmath}
\usepackage{amssymb}
\usepackage{cite}

\title[R.  Shamoyan]{On multipliers of some new analytic  function spaces in polydisc of BMOA type }

\author[1]{R. Shamoyan}

\address{Department of Mathematics, Bryansk State Technical University, Bryansk, 241050, Russia}
\email{\rm rshamoyan@gmail.com}

\date{}

\begin{document}
\maketitle  %2000 MSC numbers
 %key words and phrases

\maketitle
\begin{abstract}
We define two new BMOA type analytic function spaces in polydisk.
We provide several new results concerning coefficient multipliers of these two new BMOA analytic function spaces in polydisc.
Our results extend previously  known  assertions.
\end{abstract}
\section{Introduction}
The goal of this paper is to continue the investigation of spaces of coefficient multipliers
of some new holomorphic scales of functions in higher dimension ,namely in the unit polydisc.This note is a continuation of series of recent notes  on this topic \cite{Yaroslavcheva}, \cite {Sh1}, \cite{Sh2}, \cite{ARS}. These new spaces we study in this new note serve as  very natural extensions of the classical BMOA-type spaces in the unit polydisc. Note pointwise and coefficient multipliers of analytic spaces of Hardy and Bergman type in one and higher dimension were studied intensively by many authors during past several decades, see for example \cite{Shvedenko}, \cite{Tr1},\cite{Na1} and various references there. The investigation  of coefficient multipliers of analytic mixed norm spaces  in the unit polydisc was started in particular in recent papers of the author, see  \cite{Sh1}, \cite{Sh2}, \cite{Sh3},\cite{ARS}. Below we list standard notation and definitions which are needed and in the next section we formulate  main results of this note.To be more precise we consider two new scales of analytic functions  in the unit  polydisk of BMOA type .One of them is the dual of analytic Hardy class $H^{1}$ in polydisk and the second one is based on well-known characterization of BMOA spaces via Carleson-type measure in unit disk.Both scales in one dimension concide obviously with classical BMOA in the unit disk.Note we do not provide sharp results in this paper,
but our results are far reaching extentions of onedimensional known theorems on coefficient multipliers of BMOA type spaces.
\cite{Shvedenko},\cite{BP},\cite{JP},\cite{MP1}.  
As in many previous papers on this topic we concentrate in searching of conditions on coefficient multipliers of BMOA type space in the polydisk via
conditions on $M_{p}(f,r)$ type functions.(see definitions below).
Let $U = \{ z \in \mathbb C : |z| < 1 \}$ be the unit disc in $\mathbb C$, $\mathbb T = \partial U =
\{ z \in \mathbb C : |z| = 1 \}$, $U^n$ is the unit polydisc in $\mathbb C^n$, $\mathbb T^n \subset \partial
U^n$ is the distinguished boundary of $U^n$, $\mathbb Z_+ = \{ n \in \mathbb Z : n \geq 0 \}$, $\mathbb Z_+^n$ is the set of all multi indexes and $I = [0, 1)$.

We use the following notation: for $z = (z_1, \ldots, z_n) \in \mathbb C^n$ and $k = (k_1, \ldots, k_n) \in
\mathbb Z_+^n$ we set $z^k = z_1^{k_1} \cdots z_n^{k_n}$ and for $z = (z_1, \ldots, z_n) \in U^n$ and $\gamma \in
\mathbb R$ we set $(1-|z|)^\gamma = (1-|z_1|)^\gamma \cdots (1-|z_n|)^\gamma$ and $(1-z)^\gamma = (1-z_1)^\gamma \cdots
(1-z_n)^\gamma$. Next, for $z \in \mathbb R^n$ and $w \in \mathbb C^n$ we set $wz = (w_1z_1, \ldots, w_nz_n)$. Also, for $k \in \mathbb Z_+^n$ and $a \in \mathbb R$ we set $k + a = (k_1 + a, \ldots, k_n + a)$. For $z = (z_1, \ldots, z_n) \in \mathbb C^n$ we set $\overline z = (\overline z_1, \ldots, \overline z_n)$. For $k \in (0, +\infty)^n$ we set $\Gamma (k) = \Gamma (k_1) \cdots \Gamma (k_n)$.

The Lebesgue measure on $\mathbb C^n \cong \mathbb R^{2n}$ is denoted by $dV(z)$, normalized Lebesgue measure on $\mathbb T^n$
is denoted by $d\xi = d\xi_1 \ldots d\xi_n$ and $dR = dR_1 \ldots dR_n$ is the Lebesgue measure on $[0, +\infty)^n$.

The space of all functions holomorphic in $\mathbb U^n$ is denoted by $H(\mathbb U^n)$. Every $f \in H(\mathbb U^n)$ admits
an expansion $f(z) = \sum_{k \in \mathbb Z_+^n} a_k z^k$. For $\beta > -1$ the operator of fractional differentiation is
defined by (see \cite{DS})
\begin{equation}\label{DfFrD}
D^\beta f(z) = \sum_{k \in \mathbb Z_+^n} \frac{\Gamma(k+\beta+1)}{\Gamma(\beta+1)\Gamma(k+1)} a_k z^k,
\qquad z \in U^n.
\end{equation}

Similarly we can easily define the operator of fractional derivative for $\beta$-vector, $\beta=(\beta_1,\dots,\beta_n)$, $\beta_j>-1$, $j=1,\dots,n$.
$$D^{\vec{\beta}} f(z) = \sum_{k \in \mathbb Z_+^n} \frac{\Gamma(k_j+\beta_j+1)}{\Gamma(\beta_j+1)\Gamma(k_j+1)} a_{k_1\dots k_n} z_1^{k_1}\dots z_n^{k_n}$$

For $f \in H(\mathbb U^n)$, $0 < p < \infty$ and $r \in I^n$ we set
\begin{equation}\label{EqMpr}
M_p(f,r) = \left( \int_{\mathbb T^n} |f(r\xi)|^p d\xi \right)^{1/p},\; d\xi=d\xi_1\dots d\xi_n
\end{equation}
with the usual modification to include the case $p = \infty$. For $0<p\leq \infty$ we define standard analytic Hardy classes in the
polydisc:
\begin{equation}\label{EqDfHp}
H^p(\mathbb U^{n}) = \{ f \in H(U^n) : \| f \|_{H^p} = \sup_{r \in I^n} M_p(f,r) < \infty \}.
\end{equation}

For $n=1$ these spaces are well studied. The topic of multipliers of Hardy spaces in polydisc is relatively new, see for example \cite{BP}, \cite{Shvedenko},\cite{Tr1},\cite{Sh1},\cite{Sh2}. These spaces are Banach spaces for $p>1$ and complete metric spaces for all other positive values of $p$. Also, for $0<p\leq \infty$, $0<q<\infty$ and $\alpha > 0$ we have mixed (quasi) norm spaces, defined below as follows

\begin{equation}\label{EqMixN}
A^{p,q}_\alpha (U^n) = \left\{ f \in H(U^n) : \| f \|^q_{A^{p,q}_\alpha} = \int_{\mathbb I^n}
M_p^q(f, R) (1-R)^{\alpha q - 1} dR < \infty \right\}.
\end{equation}

If we replace the integration by $I^{n}$ above by integration by unit interval $I$ then we get other similar to these $A^{p,q}_\alpha$ analytic spaces but on subframe and we denote them by $B^{p,q}_\alpha$.
These spaces are Banach spaces for cases when both $p$ and $q$ are bigger than one, and they are complete metric spaces for all other values of parameters. We refer the reader for these classes in the unit ball and the unit disk to \cite{Shvedenko}, \cite{OF1}, \cite{OF2}, \cite{OF3} and references therein. Multipliers between $A^{p,q}_\alpha$ spaces on the unit disc were studied in
detail in \cite{Shvedenko},\cite{BP}, \cite{JP},\cite{MP1}. Multipliers of spaces on subframe are less
studied.In forthcoming paper we will study multipliers of related $B^{p,q}$ spaces.
As we see from assertions below spaces of multipliers of classes on subframe depend also from dimension $n$.

As is customary, we denote positive constants by $c$ (or $C$), sometimes we indicate dependence of a constant on a parameter by
using a subscript, for example $C_q$.

We define now one of the  main objects of this paper.

For $0<p,q<\infty$ and $\alpha > 0$ we consider Lizorkin-Triebel spaces $F^{p,q}_\alpha (U^n) = F^{p,q}_\alpha$
consisting of all $f \in H(U^n)$ such that
\begin{equation}\label{EqTrLi}
\| f \|^p_{F^{p,q}_\alpha} = \int_{\mathbb T^n} \left( \int_{I^n} |f(R\xi)|^q (1-R)^{\alpha q - 1}dR \right)^{p/q} d\xi < \infty.
\end{equation}

It is not difficult to check that those spaces are complete metric spaces, if $\min(p,q) \geq 1$ they are Banach spaces.
If we simply replace the integration by $I^{n}$ above by integration by unit interval $I$ then we get other similar to these $F^{p,q}_\alpha$ analytic spaces, but on subframe and we denote them by $T^{p,q}_\alpha$. They were studied in \cite{MitrShamoyan}.
We refer the reader to this paper for various other properties on these type classes.
We provide some new assertions on multipliers of this new $F^{p,q}_\alpha$ type spaces on subframe below. 
We note that this $F^{p,q}$ scale of spaces includes,
for $p=q$, weighted Bergman spaces $A^p_\alpha = F^{p,p}_{\frac{\alpha+1}{p}}$ (see \cite{DS} for a detailed account of these spaces). On the other hand, for $q=2$ these spaces coincide with Hardy-Sobolev spaces namely $H^{p}_\alpha =
F^{p,2}_{\frac{\alpha + 1}{2}}$, for this well known fact see \cite {Sh2}, \cite{DS} and references therein.

Finally, for $\alpha \geq 0$ and $\beta \geq 0$ we set
\begin{equation}\label{EqDinin}
A^{\infty,\infty}_{\alpha, \beta} (U^n) = \{ f \in H(U^n) : \| f \|_{A^{\infty,\infty}_{\alpha,\beta}} =
\sup_{r \in I^n} (M_\infty (D^\alpha f, r)) (1-r)^\beta < \infty \}.
\end{equation}
This is a Banach space and it is related with the well-known Bloch class in polydisc (see \cite{Shvedenko}). The Bloch class studied by many authors in various papers.
(see, for example, \cite{Shvedenko}, \cite{AP} for unit disk case and for polydisc case and also various references there).
The Bloch space is a same $$ A^{\infty,\infty}_{\alpha, \beta} (U^n)$$ space with $\alpha=\beta=1$ in definition above.
This space is also a Banach space.
For all positive values of $p$ and $s$ we introduce the following two new spaces.
We note first  replacing $q$ by $\infty$ in a usual way we will arrive at some other spaces $F^{p,\infty}$(the limit case of $F^{p,q}_\alpha$ classes.)
The limit space case $F^{p,\infty,s}(U^n)$ is defined  as a space of all analytic functions $f$ in the polydisc
such that the function $\phi(\xi) = \sup_{r \in I^n}|f(r\xi)|(1-r)^{s}$, $\xi \in \mathbb T^n$ is in $L^p(\mathbb T^n, d\xi)$, $s\geq 0$, $0<p\leq\infty$.

Finally, the limit space case $A^{p,\infty,s}(U^n)$ is the space of all analytic functions $f$ in the polydisc such that (see \cite{Shvedenko})
$$\sup_{r\in I^n}M_{p}(f,r)(1-r)^{s}<\infty, s\geq 0, 0<p\leq\infty$$

Obviously the limit $F^{p,\infty,s}$ space is embedded in the last space we defined.  These last two spaces are Banach spaces for all $p\geq 1$ and they are complete metric spaces for all other positive values of $p$.

The following definition of coefficient multipliers is well known in the unit disk. We provide a natural extension to the
polydisc setup.This definition can be seen in  papers of the first author before .
\begin{definition}Let $X$ and $Y$ be  quasi normed subspaces of $H(U^n)$. A sequence $c = \{ c_k \}_{k \in \mathbb Z_+^n}$ is said to be a coefficient multiplier from $X$ to $Y$ if for any function $f(z) = \sum_{k \in \mathbb Z_+^n} a_k z^k$ in $X$ the function $h = M_cf$ defined by $h(z) = \sum_{k \in \mathbb Z_+^n} c_k a_k z^k$ is in $Y$. The set of all multipliers from $X$ to $Y$ is denoted by $M_T(X, Y)$.
\end{definition}

The problem of characterizing the space of multipliers (pointwise multipliers and coefficient multipliers)
between various spaces of analytic functions is a classical problem in complex function theory,
there is vast literature on this subject, see \cite{BP}, \cite{DS}, \cite{JP}, \cite{MP1}, and references therein.

In this paper we are looking for some extensions of already known classical theorems, namely we are interested in spaces of multipliers acting into analytic BMOA type spaces in the unit polydisc and from these spaces into certain well studied classes like mixed norm spaces, Bergman spaces and Hardy spaces. We note that the analogue of this problem of description of spaces of multipliers in $\mathbb R^n$ for various functional spaces in $\mathbb R^n$ was considered previously by various authors in recent decades (see \cite{Shvedenko}).
The intention of this note is to study the spaces of coefficient multipliers of new analytic function spaces in polydisc.
As we indicated this topic is well-known and various results by many authors were published in recent decades in this area starting from
classical papers of Hardy- and Littlewood (see for example \cite{Shvedenko}) for many known results  and various references there.
Nevertheless the study of coefficient multipliers of analytic function spaces in higher dimension namely in polydisc is a new area of research and only several sharp results are known till now  see ,for example, \cite{Shvedenko},\cite{Yaroslavcheva},\cite{Sh1} -- \cite{Tr1}  and various references there. We complement these results and our results from \cite{Sh1},\cite{Sh2},\cite{Sh3} in this note providing new theorems on coefficient multipliers of some new analytic function spaces in polydisc.

We define also by  $dm_2(\xi)$ the Lebesgue measure on the unit disk and replacing $2$ by $2n$ the Lebesgue measure on the polydisc. 

Let $L^{\vec{p}}(\vec{\alpha})$ be the space of all measurable functions in $U^n$ so that
\[\|f\|_{L^{\vec{p}}(\vec{\alpha})}=\Biggr(\int\limits_{U}(1-|\xi_n|)^{\alpha_n}\Biggr(\int\limits_{U}(1-|\xi_{n-1}|)^{\alpha_{n-1}}...\Biggr(\int\limits_{U}(1-|\xi_{1}|)^{\alpha_{1}}|f(\xi_1, ..., \xi_n)|^{p_1}\times\]\[\times dm_2(\xi_1)\Biggl)^{\frac{p_2}{p_1}}...dm_2(\xi_{n-1})\Biggl)^{\frac{p_n}{p_{n-1}}}dm_2(\xi_{n})\Biggl)^{\frac{1}{p_{n}}}<+\infty.\]
where all $p_j$ are  positive and all $\alpha_j>-1$ for any $j=1,...,n$
Let $H^{\overrightarrow{p}}(\overrightarrow{\alpha})=L^{\vec{p}}(\vec{\alpha})\cap H(U^n)$. These spaces are one of the main objects of study of this paper (see \cite{Yaroslavcheva}, \cite{YaroslavchevaShamoyan} for properties of these new classes).

We define now new mixed norm weighted Hardy analytic spaces in polydisc in the following natural way.

Let $$H^{\overrightarrow{p}}_{\overrightarrow{\alpha}}=\sup\limits_{r_j<1}\Biggr(\int\limits_T...\Biggr(\int\limits_T\Biggr(\int\limits_T|f(r\overrightarrow{\xi})|^{p_1}d\xi_1\Biggl)^{\frac{p_2}{p_1}}d\xi_2\Biggl)^{\frac{p_3}{p_2}}\dots d\xi_n\Biggl)^{\frac{1}{p_n}}\prod\limits_{j=1}^{n}(1-r_j)^{\alpha_j};$$ $p_j\in(0; \infty]$; $\alpha_j\geq0$, $j=1,\dots,n$.
Note in unit disk we get well-studied classical weighted Hardy spaces \cite{Shvedenko}, \cite{DS}, \cite{GL}, \cite{JP}, \cite{ARS}, $H_0^{p,\dots,p}=H^p(U^n)$.
The main objects of this paper are two new BMOA type spaces in higher dimension.
BMOA type spaces in higher dimension can be defined differently
starting from their various onedimensional characteristics .
\cite{Shvedenko} and references there.
We first define BMOA in polydisk as dual space of classical analytic hardy space $H^{1}$
in polydisk which we mentioned above.We note this serves as direct extention of BMOA spaces from the unit disk 
\cite{Shvedenko}.Another natural way to extend directly
BMOA spaces from unit disk to higher dimension is  to use their characterizations in unit disk via Carleson type measures,this type approach was used already by first author before.Namely
we define the other main object of this paper as follows.
Let $F^{q,\infty}_{s,\alpha}$ be a space of analytic
functions in the unit polydisk ,so that 

$$\int_{U^{n}}\frac{|f(z)|^{q}(1-|z|)^{s}(1-|w|)^{\alpha}}{|1-wz|^{\alpha+1}}dm_{2n}(z)$$
 is a bounded function in the unit polydisk as a function of $w$ variable.
We assume that $s>-1$ and other values of parameters above are positive numbers.
Note both spaces coincide in the unit disk with each other and with classical BMOA space.\cite{Shvedenko}.These spaces are the main objects of this paper.

\section{Preliminary results}

For formulation of all our main results of this note we need several lemmas, almost all of them are known and taken from recently published papers. To prove our results we need to use a 
simple idea.We connect our BMOA type analytic function spaces in polydisk with much simpler well- studied spaces in the polydisk and using these connections(embeddings)we find information about sizes of spaces of coefficient multipliers we search for.
Another important ingredients of our work are estimates of Bergman kernel in polydisc in some new spaces with usual and unusual quazinorms in polydisc.Note much information concerning this was taken from recent papers of the first author on this topic.see \cite{MitrShamoyan},\cite{KurShamoyan} .As usual we define positive constants by $C$ (or $c$) with various (mainly lower) indexes depending on situation.

In particular  we need the following lemma 1:
\begin{lemma}[see \cite{Yaroslavcheva}, \cite{YaroslavchevaShamoyan}]
Let $g(z_1, ...,z_n)=\prod\limits_{j=1}^{n}\frac{1}{(1-z_j)^{\beta_j}}, g_r(z)=g(rz), r\in I^n, z_j \in U, \beta_j>0, j=1, ...,n$, $g_r(z)=g(rz)$. Then for some large enough $\beta_0$ and all $\beta_j>\beta_0, j=1, ...,n$ we have
\[\|g_r\|_{H^{\overrightarrow{p}}(\overrightarrow{\alpha})}\leq c\prod\limits_{j=1}^{n}(1-r_j)^{\frac{\alpha_j}{p_j}-\frac{\beta_jp_j-2}{p_j}},\]
for $\alpha_j>-1$ and all $0<p_j<\infty, r_j \in I, j=1, ..., n;$
for $p_j=\infty$ we have $\beta_j>0, j=1, ..., n$ and
\[\|g_r\|_{H^\infty}\leq C \prod\limits_{j=1}^{n}(1-r_j)^{-\beta_j}, \ \ r_j\in (0;1).\]
for all $j=1,...,n$
\end{lemma}

One of the main goals of this note is to study multipliers of some  new H-type analytic so-called mixed norm spaces with unusual norms and quazinorms. For that reason first we find direct connections of these new
analytic spaces with analytic function spaces with simpler norms and quazinorms based on various simple embeddings and then use this connection and  already known information about coefficient multipliers of analytic classes with simpler quazinorms to get information we search for in standard terms of Ñlassical $M_{p}(f,r)$ function. This shows the importance of embeddings connecting various analytic function spaces in polydisc.We in particular need these assertions.
\begin{lemma}[see \cite{Yaroslavcheva}, \cite{YaroslavchevaShamoyan}]
Let $f\in H^{\overrightarrow{p}}(\overrightarrow{\alpha}), 0<p_j\leq1, j=1, ..,n, \alpha_j>-1, \gamma_j=2-p_j, j=1, ..., n.$. Then we have
\[\int\limits_0^1(1-r_n)^{\alpha_{n}+\gamma_{n}-1}...\Biggr(\int\limits_0^1(1-r_2)^{\alpha_{2}+\gamma_{2}-1}\int\limits_0^1(1-r_1)^{\alpha_{1}+\gamma_{1}-1}\times\]\[\times M^{p_1}_1(r_1, ..., r_n, f)r_1 dr_1\Biggl)^{\frac{p_2}{p_1}}r_2dr_2\Biggl)^{\frac{p_3}{p_2}} ...r_ndr_n\leq c\|f\|_{H^{\overrightarrow{p}}(\overrightarrow{\alpha})}^{p_n};\]
where \[M^{p}_p(r_1, ..., r_n, f)=\int\limits_{T^n}|f(r_1\xi_1, ..., r_n\xi_n)|^pd\xi_1 ...d\xi_n;\]
and for all $0<p_j<\infty$, $\alpha_j>-1$, $j=1,\dots, n$
\[\sup\limits_{z\in U^n}|f(z)|\prod\limits_{j=1}^{n}(1-|z_{j}|)^{\frac{2}{p_j}+\frac{\alpha}{p_j}}\leq c\|f\|_{H^{\overrightarrow{p}}(\overrightarrow{\alpha})}.\]
\end{lemma}

From lemma 2 we get
\begin{corollary}
 Let $f\in H^{\overrightarrow{p}}({\overrightarrow{\alpha})}, 0<p_j\leq1, j=1, ...,n.$ Then $M_1(f, r)(1-r)^\tau<\infty$, $r\in I$, where $\tau=\sum\limits_{j=1}^{n}(\alpha_j-p_j+2)\frac{1}{p_j}, \alpha_j>p_j-2, j=1, ...,n.$
\end{corollary}

Similar result is true for $M_1(f,r)$, $r\in I^n$.Here we just must replace $M_1(f,r)$ by $M_1(f,r_{1},\ldots,r_{n})$ and $(1-r)^{\tau}$ by $\prod(1-r_{j})^{\tau_j}$ ,where 
$\tau_{j}p_{j}=\alpha_j-p_j+2$, and where $r_{j}\in I$.
We will use this corollary and this remark constantly later.
We introduce now new spaces of analytic functions in polydisc and provide estimates of Bergman kernel in quazinorms of these spaces in polydisk.
Let for all positive values of indexes $p,q,\gamma$
\[M^{p, q}_{\gamma}=\{f \in H(U^n): \|f\|^p_{M^{p,q}_\gamma}=\int\limits_T\Bigr(\int\limits_0^1\int\limits_0^1|f(r\zeta)|^q(1-r)^{\gamma q-1}dr_1 ...dr_n\Bigl)^{\frac{p}{q}}d\zeta <\infty\}\]
\[ 0<p,q<\infty,\; \gamma>0.\]
These spaces were studied recently in \cite{MitrShamoyan}.

Note simply replacing $T$ by $T^n$ we get immediately known analytic $F^{p,q}_\alpha$ spaces (see \cite{ARS}).
We define also $\|f\|_{M^{p, \infty}_{\gamma}}=\Biggr(\int\limits_T\Bigr(\sup\limits_{r\in I^n}|f(r\zeta)|(1-r)^{\alpha}\Bigl)^pd\zeta\Biggl)^{\frac{1}{p}}<\infty$, $0<p<\infty$, $\alpha\geq 0$ and with usual modification for $p=\infty, 0<q<\infty$. We put $A/p=0$ if $p=\infty$.
In the following several lemmas we do not specify parameters
$p,q$ separately meaning they can accept all positive  values including 
$\infty$.We also do not put any restriction on parametrs
$\gamma$ and $s$ except those which are accepted in definition
of spaces which they represent in following lemmas.
\begin{lemma}[see \cite{MitrShamoyan}]
Let
\[g_R(z)=\frac{1}{(1-R z)^{\beta+1}}, {R}\in I, {z}\in U^{n}.\]
Then
\[ \|g_{R}\|_{M^{p, q}_{\gamma}}\leq \frac{C_1}{(1-{R})^{n(\beta+1)-\gamma n-\frac{1}{p}}},  \beta>\gamma+\frac{1}{np}-1,\]

\[ \|g_{R}\|_{M^{p,\infty}_{\gamma}}\leq \frac{C_2}{(1-{R})^{n(\beta+1)-\gamma n-\frac{1}{p}}},  \beta>\gamma+\frac{1}{np}-1, \]

\[ \|g_{R}\|_{M^{ \infty,q}_{\gamma}}\leq \frac{C_3}{(1-{R})^{({\beta}+1-{\gamma})n}},  \beta>\gamma-1.\]

\end{lemma}
We note Lemma 3 in combination with standard arguments based on so-called closed graph theorem  allows to show immediately the Theorem 4 below with the help of lemma 9 and 10.
Now we formulate a complete analogue of lemma 3 for $T^{p, q}_\gamma$ spaces in polydisc.
\begin{lemma}[see \cite{MitrShamoyan}]
Let
\[g_R(z)=\frac{1}{(1-R z)^{\beta+1}}, {R}\in I^n, {z}\in U^{n}, \beta>0.\]
Then
\[ \|g_{R}\|_{T^{p, q}_{\gamma}}\leq \frac{C_1}{(1-{R})^{\beta+1-\frac{\gamma}{n}-\frac{1}{p}}}, \beta>\frac{\gamma}{n}+\frac{1}{p}-1,\]

\[ \|g_{R}\|_{T^{\infty, q}_{\gamma}}\leq \frac{C_2}{(1-{R})^{\beta+1-\frac{\gamma}{n}}},  \beta>\frac{\gamma}{n}-1,\]

\[ \|g_{R}\|_{T^{p, \infty}_{\gamma}}\leq \frac{C_3}{(1-{R})^{\beta+1-\frac{\gamma}{n}-\frac{1}{p}}},  \beta>\frac{\gamma}{n}+\frac{1}{p}-1.\]
\end{lemma}

The proof of lemma 5 is easy.We omit it leaving it to readers.

\begin{lemma}
For $g_R(z)=\prod\limits_{j=1}^{n}\frac{1}{(1-R_jz_j)^{\beta_j+1}}, R\in I^n, z\in U^n,$
we have
\[\|g_R\|_{H^{\overrightarrow{p}}_{\overrightarrow{\alpha}}}\leq\frac{C}{\prod\limits_{j=1}^{n}(1-R_j)^{\beta_j+1-\frac{1}{p_j}-\alpha_j}}, 
\beta_j>\alpha_j+\frac{1}{p_j}-1, j=1, ..., n.\]
\end{lemma}

The following is a complete analogue of lemma 3 for  $A^{p, q}_\gamma, F^{p, q}_\gamma, H^s$ spaces,
namely we provide estimates for Bergman kernel in polydisc in these spaces .

\begin{lemma}[see \cite{ARS}]
Let
\[g_R(z)=\frac{1}{(1-R z)^{\beta+1}}, {R}\in I^n, {z}\in U^{n}, \beta>0.\]
Then
\[ \|g_{R}\|_{A^{p, q}_{\gamma}}\leq \frac{C_1}{(1-{R})^{\beta-{\gamma}-\frac{1}{p}+1}}, \beta>\gamma+\frac{1}{p}-1,\]

\[ \|g_{R}\|_{F^{p, q}_{\gamma}}\leq \frac{C_2}{(1-{R})^{\beta-{\gamma}-\frac{1}{p}+1}}, \beta>\gamma+\frac{1}{p}-1,\]

\[ \|g_{R}\|_{H^s}\leq \frac{C_3}{(1-{R})^{\beta+1-\frac{1}{s}}},  \beta>\frac{1}{s}-1.\]
\end{lemma}
The following is a complete analogue of lemma 3 for  $A^{p, \infty}_\gamma, A^{\infty, q}_\gamma, F^{p, \infty}_\gamma, H^\infty, Bl$
namely we provide estimates for Bergman kernel in polydisc in these spaces .

\begin{lemma}[see \cite{ARS}]
Let
\[g_R(z)=\frac{1}{(1-R z)^{\beta+1}}, {R}\in I^n, {z}\in U^{n}, \beta>0.\]
Then
\[ \|g_{R}\|_{A^{p, \infty}_{\gamma}}\leq \frac{C_1}{(1-{R})^{\beta-{\gamma}-\frac{1}{p}+1}}, \beta>\gamma+\frac{1}{p}-1,\]

\[ \|g_{R}\|_{A^{\infty, q}_{\gamma}}\leq \frac{C_2}{(1-{R})^{\beta-{\gamma}+1}}, \beta>\gamma-1,\]

\[ \|g_{R}\|_{F^{p, \infty}_{\gamma}}\leq \frac{C_3}{(1-{R})^{\beta-{\gamma}-\frac{1}{p}+1}}, \beta>\gamma+\frac{1}{p}-1,\]

\[ \|g_{R}\|_{Bl}\leq \frac{C_4}{(1-{R})^{\beta+1}}, \beta>-1,\]

\[ \|g_{R}\|_{H^\infty}\leq \frac{C_5}{(1-{R})^{\beta+1}}, \beta>-1.\]
\end{lemma}

Next lemma plays the same role as lemma 2 for $H^{\overrightarrow{p}}({\overrightarrow{\alpha}})$.
\begin{lemma}[see \cite{ShamoyanMitrich}]
Let $0<p_j<1, j=1,...,n.$ Then
\[\int\limits_{U^n}|f(z_1, ..., z_n|\prod_{j=1}^{n}(1-|z_j|)^{\frac{1}{p_j}-2}dm_{2n}(z)\leq c \|f\|_{H^{\overrightarrow{p}}}=\]\[=\sup\limits_{0<r_j\leq1}
\int\limits_T\Biggr(\int\limits_T ... \Biggr(\int\limits_T|f(r_1\zeta_1, ..., r_n\zeta_n)|^{p_1}d(\zeta_1)\Biggl)^{\frac{p_2}{p_1}}...d(\zeta_n)\Biggl)^{\frac{1}{p_n}}.\]
\end{lemma}

We get the following assertion from lemma 8 immediately using standard arguments.

\begin{corollary}
Let $p_j\in (0,1)$, $j=1,\dots,n$, then
\[\sup_{r\in I^n}M_1(f,r)\prod_{k=1}^{n}(1-r_k)^{\frac{1}{p_k}-1}\leq c\|f\|_{H^{\vec{p}}}\]
\end{corollary}
\begin{corollary}
Let $p_j\in (0,1)$,$\alpha_{j}>0$, $j=1,\dots,n$, then
\[\sup_{r\in I^n}M_1(f,r)\prod_{k=1}^{n}(1-r_k)^{\frac{1}{p_k}-1+\alpha_{k}}\leq c\|f\|_{H^{\vec{p}}_{\alpha}}\]
\end{corollary}

Weighted Hardy space case $H^{\vec{p}}_{\vec{\alpha}}$ follows from corollary 2 if we apply it to $f_R, R\in I$ function using then standard arguments.
We add here two important remarks.
As we will see this corollary 2 and it is weighted version we just mentioned can be used immediately instead of corollary 1 to get same theorems  for $H^{\vec{p}}_{\vec{\alpha}}$ instead of $H^{\vec{p}}(\vec{\alpha})$ for those cases where $M_T(Y,H^{\vec{p}}(\vec{\alpha}))$ appears  in discussions after theorems below, we leave this procedures to readers. (we denote here   by $Y$ one of BMOA type spaces we study in this paper)
 On the other hand Lemma 8 gives ways  to find easily simple necessary conditions for $g$ function in usual in these issues terms of $M_1(f,r)$ function to be a multiplier from $Y$ to $H^{\overrightarrow{p}}_{\vec{\alpha}}$ for BMOA type $Y$ spaces we considered above (also with the help of lemma's 3-7 and the closed graph theorem).We will also see this below easily.
Our last two assertions are the most important lemmas for this work.They connect BMOA type spaces with much simpler analytic function spaces in unit polydisk which leads to solutions of some problems we already indicated for this work.
\begin{lemma}
Let $\alpha>0$,$s>0$,$0<q<\infty$ 
Then if $f\in A^{\infty,q}_{s}$ 
then $$f\in F^{\infty,q}_{sq-1,\alpha}$$ with related
estimates for quazinorms.
and if $f\in F^{\infty,q}_{sq-1,\alpha}$ then  $$f\in A^{\infty,\infty}_{s}$$ with related estimates for quazinorms.
\end{lemma}
\begin{lemma}
Let $f\in BMOA$ then $$f \in A^{\infty,\infty}_{\beta,\beta}$$ for all $\beta>0$.
If $$f \in A^{\infty,\infty}_{\beta,\beta-1-\alpha}$$ for all $\beta>\alpha+1$ and $\alpha>-1$ then $f\in BMOA$
with related estimates for quazinorms.
\end{lemma}
The $F^{\infty,q}_{s,\alpha}$ analytic spaces can also be defined in a more general for  vectors $s$ and $\alpha$
when parameters related with weights are vectors.
Most  results  below  concerning  this $F^{\infty,q}_{s,\alpha}$  space  can  be  reformulated  in  this  more  general  form.
We  add  comments concerning  proofs  of  these  assertions.
The first part of lemma 9 is obvious.For unit disk case it is 
known ,see ,for example, \cite{Shvedenko}.The second part of lemma 9 for the case of unit ball can be seen with a complete proof in \cite{OF1}.
The base of that proof is a Bergman reproducing formula and 
Forelli-Rudin type estimates.Note the polydisk version of that lemma is also valid and the proof is a copy of a proof of unit ball case.We omit details referying the reader to \cite{OF1}
for all details.Turning to the proof of lemma 10 we note the proof of the first part can be seen in \cite{Shamoyan1}
and the proof of the second one follows from the following
simple argument.Note first for any $\alpha$, $\alpha>0$ we have 
$f \in A^{1,1}_{\alpha}$ for any $f$ 
function so that $f\in H^{1}$.Since both classes are Banach spaces this means that the reverse inclusion holds for dual spaces.But the dual of Bergman class  is well-known in polydisk
see for example \cite{Shvedenko} and hence the assertion of lemma follows.
Using last lemmas estimates of Bergman kernel
and $g_R$ function (see also similar assertions above) can be found for
BMOA- type spaces we introduced in this paper.
These are estimates of $g_{R}$ for 
$A^{\infty,\infty}_{\alpha}$
or $A^{\infty,q}_{\alpha}$ spaces which can already be seen in this section.
Note now all assertions of next section concerning two scales of BMOA type functions we study in this paper  are in some sence parallel to each other in a sense that each assertion on BMOA as we defined above leads to similar assertion concerning 
$F^{q,\infty}_{s,\alpha}$ analytic spaces with similar arguments for proofs.We just have to replace embeddings in last two lemmas.This remark can be applied to each assertion
below formulated for $F^{q,\infty}_{s,\alpha}$  spaces in polydisk.  
The exposition in this note is sketchy.We alert the reader in advance
many results of this note( and even proofs which are easy calculations) are left to interested readers, but all details of schemes of proofs  for that are indicated below by us.  
We introduce a vital notation for this paper.

\section{ Main Results on Coefficient Multipliers in some new spaces of analytic functions in polydisc of BMOA - type}

The intention of this section to provide main results of this paper.The base of all proofs is a standard argument based on closed graph theorem
used before by many authors (see for example \cite{Shvedenko}) in one dimension and in higher dimension (see for example \cite{Sh1},\cite{Sh2},\cite{Sh3}, \cite{ARS} and references there). This in combination with various lemmas  containing estimates of Bergman kernel in some new unusual analytic spaces we provided above leads directly to completion of proofs. Note again this scheme in $U^n$ was used by first author also in previous mentioned papers
\cite{Sh1},\cite{Sh2},\cite{Sh3}, \cite{ARS}. The appropriate definitions of coefficient multipliers in polydisc and Bergman kernels in polydisc serve as important steps for solutions.
Our first theorem and the second theorem provide new necessary conditions for multipliers into various known mixed norm spaces from BMOA-type spaces which we study in this paper based on  first lemmas. Note these two theorems are not new  in case of unit disk (see \cite{Shvedenko}). So they can be viewed as far reaching extensions of these old results on coefficient multipliers of  BMOA-type(see \cite{Shvedenko})in the unit disk.
The following theorem can be seen also as a corollary of remark  above which provides direct ways to get estimates of Bergman kernel in mentioned BMOA-type spaces and estimates of Bergman kernel for $A^{p,q}_\alpha$ and $F^{p,q}_\alpha$ spaces. We alert the reader many arguments  in this paper are sketchy  ,this allows using short space to provide
big amount of new results on multipliers of BMOA type spaces in polydisk.We indicate also for that reason only shortly changes that are needed in proofs we provide below to get new assertions.

We  formulate  first  two  assertions  for  $F^{\infty,q}_{sq-1,\alpha}$   spaces  in  this  paper .

\begin{theorem}
1) Let $0<p,v<\infty$,$\gamma>0$, $g\in H(U^n)$.Let also $g(z)=\sum\limits_{|k|\geq0}^{\infty}c_{k_1, ..., k_n}z_1^{k_1}...z_n^{k_n}$, be a multiplier 
from $F^{\infty,q}_{sq-1,\alpha}$ to $A^{p, v}_\gamma(U^n)$ then
\[\sup\limits_{r\in I^n}M_p(D^\beta g, r)\prod\limits_{j=1}^{n}(1-r_j)^{\gamma+\beta_j+1-s}<\infty,\]
for all $\beta_j > s-1-\gamma$, $j=1, ..., n$.

2)Let $0<p,v< \infty$,$\gamma>0$, $g\in H(U^n)$.Let also $g(z)=\sum\limits_{k_1\geq0, ..., k_n\geq0}^{\infty}c_{k_1, ..., k_n}z_1^{k_1}...z_n^{k_n}$ be a multiplier from 
$F^{\infty,q}_{sq-1,\alpha}$ to $F^{p, v}_\gamma(U^n)$ then for $v\leq p$
\[\sup\limits_{r\in I^n}M_p(D^\beta g, r)\prod\limits_{j=1}^{n}(1-r_j)^{\gamma+\beta_j+1-s}<\infty,\]
for all $\beta_j>s-1-\gamma$,$ j=1, ..., n$ .
\end{theorem}

\begin{theorem}
1) Let $0<p,q<\infty$,$\gamma>0$, $g\in H(U^n)$.Let also $g(z)=\sum\limits_{|k|\geq0}^{\infty}c_{k_1, ..., k_n}z_1^{k_1}...z_n^{k_n}$, be a multiplier from $BMOA$ to $A^{p, q}_\gamma(U^n)$ then
\[\sup\limits_{r\in I^n}M_p(D^\beta g, r)\prod\limits_{j=1}^{n}(1-r_j)^{\gamma+\beta_j+2+v}<\infty,\],where $\beta>-v-2$
for any,$v$, $v>-1$ .

2)Let $0<p,q<\infty$, $g\in H(U^n)$, $g(z)=\sum\limits_{k_1\geq0, ..., k_n\geq0}^{\infty}c_{k_1, ..., k_n}z_1^{k_1}...z_n^{k_n}$ be a multiplier from $BMOA$ to $F^{p, q}_\gamma(U^n)$ then for $q\leq p$
\[\sup\limits_{r\in I^n}M_p(D^\beta g, r)\prod\limits_{j=1}^{n}(1-r_j)^{\gamma+\beta_j+2+v}<\infty,\],where $\beta>-v-2$
for any $v$, $v>-1$.
\end{theorem}
The following theorem also can be seen  as a corollary of remarks after lemma 9 and 10 above which provides estimates of Bergman kernel in  mentioned BMOA- type  spaces and the closed graph theorem.
Let below Y be again the $F^{\infty,q}_{sq-1,\alpha}$ space which we introduced above.
The case of BMOA needs small changes which we already indicated.
\begin{theorem}
Let $p,\gamma>0$ .
Let $g\in H(U^n)$, $g(z)=\sum\limits_{|k|\geq 0}c_k z^k$
1)If $g \in M_T(Y; F^{p, \infty}_\gamma)$ then
\[\sup\limits_{r\in I^n}M_p(D^\beta g, r)\prod\limits_{j=1}^{n}(1-r_j)^{\gamma+\beta_j+1-s}<\infty,\]
for all $\beta_j > s-1-\gamma$, $j=1, ..., n$ .

2)If $g \in M_T(Y; A^{p, \infty}_\gamma)$ then
\[\sup\limits_{r\in I^n}M_p(D^\beta g, r)\prod\limits_{j=1}^{n}(1-r_j)^{\gamma+\beta_j+1-s}<\infty,\]
for all $\beta_j > s-1-\gamma$ $j=1, ..., n$.

\end{theorem}
\begin{remark}
Similarly conditions on $M_T(Y,X)$ for $X=T^{p,q}_\gamma$ and $B^{p,q}_\gamma$ can be found for $p<q$ or $p=q$. We have to use in addition that $T^{p,q}_\gamma\subset B^{p,q}_\gamma$, $p\leq q$ and the fact that
\[\sup_{r\in I} M_p(f,r)(1-r)^\gamma\leq c\|f\|_{B^{p,q}_\gamma},\; \gamma>0,\; 0<p,q\leq \infty.\]
These estimates are easy to see or can be seen in various
textbooks on Bergman spaces.
\end{remark}

Note that we can also replace  $T^{p,q}_{\alpha}$ spaces in previous remark we just clarified by 
$H^{\vec{p}}_{\vec{\alpha}}$ classes .To formulate   this new assertion  
we have to replace one estimate on $T^{p,q}_{\gamma}$
spaces we just mentioned by another similar estimate from below via $M_{p}(f,r)$ function for $H^{\vec{p}}_{\vec{\alpha}}$ classes .Both can be seen in previous section. For $H^{\overrightarrow{p}}(\overrightarrow{\alpha})$  spaces  the situation is the same and similar arguments can be easily applied.
In both cases we need estimates from below via $M_{p}(f,r)$
function.
Note changing an embedding we use in proofs of theorems 1,2 for BMOA by
similar embedding for $F^{\infty,q}_{s,\alpha}$ (or reverse lemma 10 by lemma 9) we can easily  formulate assertions of theorems 1 and 2 directly in terms of
for $F^{\infty,q}_{s,\alpha}$ analytic classes or by BMOA in the unit polydisk.We leave this last easy procedure to interested readers.

Proofs of theorems 1 and 2.

Here is a complete scheme of proof of various assertions of type contained in theorem 1 and theorem 2.

Assume $\{c_k\}_{k\in Z^n_+}\in M_T(Y, X)$ where $$X=A^{p,q}_\alpha \;or\; F^{p,q}_\alpha, 0<p,q\leq \infty, \alpha>0$$ or
$$X=T^{p,q}_\alpha \;or\; M^{p,q}_\alpha, 0<p,q\leq \infty, \alpha>0,$$
or $X=Bl,\; H^\infty$.

An application of the closed graph theorem gives $\|M_cf\|_X\leq c\|f\|_{Y}$.

Let $w\in U^n$ and set $f_w(z)=\frac{1}{1-wz}$, $g_w=M_cf_w$. Then we have 
$$D^\beta(g_w)=D^\beta M_cf_w=M_c(D^\beta f_w)=c\left(M_c\frac{1}{(1-wz)^{\beta+1}}\right);\; \beta>0$$
which gives
$$\|D^\beta g_w\|_X\leq c\|\tilde{f}\|_{Y}=c\left\|\frac{1}{(1-wz)^{\beta+1}}\right\|_{Y}.$$

This give immediately by lemma above an estimate
$$\|D^\beta g_w\|_X\leq c\frac{1}{\prod\limits_{j=1}^n(1-r_j)-\tau_{0}}.$$ ,for some fixed $\tau_{0}$
depending on parametrs involved. 
We finish easily the proof using only one known embedding (depending each time on quazinorm of $X$ space).

For $X=A^{p,q}_\alpha$ we use embedding 
$$\sup_{r\in I^n}(M_p(f,r)(1-r)^\alpha\leq c\|f\|_{A^{p,q}_\alpha},\; f\in H(U^n),\; 0<p\leq\infty,\;\alpha\geq 0.$$

For $X=F^{p,q}_\alpha$, $p\leq q$ we use the same embedding with inclusion 
$$F^{p,q}_\alpha\subset A^{p,q}_\alpha, \; p\leq q.$$
For $p=\infty$ or $q=\infty$ we have to replace above $p$ by $\infty$ or $q$ by $\infty$. The case if $X=Bl$ or $X=H^\infty$ we have to modify the scheme we provided above.

For $X=T^{p,q}_\alpha$ we use same arguments and the embedding
$$\sup_{r\in I}M_p(f,r)(1-r)^\alpha\leq c\|f\|_{T_\alpha^{p,q}};\; 0<p,q\leq\infty;\; p\leq q;\; \alpha>0,\; f\in H(U^n)$$
with obvious modification for $p=\infty$ or $q=\infty$.
Note that the same embedding is true for $B^{p,q}_{\alpha}$

As another corollary of Lemma's above we have the following theorems 3 and 4, 5. They provide various necessary conditions for a function to be a multiplier from mixed norm well-studied analytic Besov type and new Lizorkin-Triebel type spaces and
$H$-type new Hardy- type and Bergman- type analytic spaces with mixed norm into new BMOA type spaces in polydisk we defined and study in this paper namely BMOA and $F^{\infty,q}_{\alpha,s}$ analytic function spaces.These are far reaching extentions of previously known results.see \cite{Shvedenko}
The crucial embeddings for BMOA type spaces we indicated for these two spaces at the end of our previous section again as previously in theorems 1,2 serve as base of proofs of all our results we will formulate below.(theorems 3,4,5.)
Note as in previous case (theorems 1 and 2 )we consider  
analytic BMOA type spaces simultianiously , since the case of $F^{q,\infty}_{\alpha,s}$ and BMOA
spaces are similar and proofs are parallel, we have to change
only an embedding for BMOA by embedding for  
$F^{q,\infty}_{\alpha,s}$ which can be seen in lemmas 9 and 10 above.
Note some results contained in theorems below  are new even in case of unit disk, on the other hand  others can be viewed as direct extensions of well-known one dimensional results to the higher dimensional case of polydisc and we also record them. Note we  here  as above include in some cases "limit" spaces when some parameters involved are equal to infinity. We also leave some cases open for interested readers since they do not contain any new ideas for proofs.
The following theorem can be seen  as a corollary of lemma's above  which provides estimates of Bergman kernel in  mentioned mixed norm   spaces $A^{p, q}_\gamma$,$F^{p, q}_\gamma$ and the closed graph theorem and lemma 9 and 10
Various concrete values of $\beta_{0}$ depending on various parameters and on a concrete space can be also easily recovered by readers in all theorems of this note.
We also leave this to readers.
We formulate below three theorems for 
$F^{\infty,q_{1}}_{vq_{1}-1,\alpha}$ spaces $q_{1}\in(0,\infty)$,$v>0$,
$\alpha>0$.We consider a little bit more general case of just mentioned spaces when
$v$ is a vector from components $v_{j}$.(see for this a remark above)
We also will denote it by Y below to avoid big amount of indexes
in formulations of theorems.
Note $q_{1}$ and $\alpha$ do not even appear below in spaces of coefficient multipliers we found .
Note also the base of proof is lemma 9 and the closed graph theorem.Similar assertions (complete analogues of theorems below ) are valid also for analytic BMOA space we defined above.We have to replace only lemma 9 by lemma 10 using same arguments in our proofs.
We leave this procedure  to interested readers.
We assume always below that $\frac{a}{b}=0$,if $b=\infty$.

\begin{theorem}
Let $p>0$,$q>0$,$\gamma>0$.
and let $g(z)=\sum\limits_{k \in \mathbb Z_+^n}c_kz^k, g\in H(U^n),$ then 
1)If $g\in M_T(A^{p, q}_\gamma, Y)$ then
\[\sup\limits_{r\in I^n}M_\infty(D^{\beta}g, r_1, ..., r_n)\prod\limits_{j=1}^{n}(1-r_j)^{\tau_j}<\infty, \] $ \tau_{j}=v_{j}+\beta-\gamma-\frac{1}{p}+1$,$\beta>\max\limits_{j}\left(\gamma+\frac{1}{p}-1-v_{j}\right).$

2)If $g\in M_T(F^{p, q}_\gamma, Y)$ then
\[\sup\limits_{r\in I^n}M_\infty(D^{\beta}g, r_1, ..., r_n)\prod\limits_{j=1}^{n}(1-r_j)^{\tau_j}<\infty, \] $ \tau_j=v_{j}+\beta-\gamma-\frac{1}{p}+1, j=1, ..., n, \beta>\max\limits_{j}\left(\gamma+\frac{1}{p}-1-v_{j}\right).$

3)If $g\in M_T(F^{p, \infty}_\gamma, Y)$ then
\[\sup\limits_{r\in I^n}M_\infty(D^{\beta}g, r_1, ..., r_n)
\prod\limits_{j=1}^{n}(1-r_j)^{\tau_j}<\infty, \] $ \tau_j=v_{j}+\beta-\gamma-\frac{1}{p}+1, j=1, ..., n, \beta>\max\limits_{j}\left(\gamma+\frac{1}{p}-1-v_{j}\right).$

4)If $g\in M_T(A^{p, \infty}_\gamma, Y)$ then
\[\sup\limits_{r\in I^n}M_\infty(D^{\beta}g, r_1, ..., r_n)
\prod\limits_{j=1}^{n}(1-r_j)^{\tau_j}<\infty, \] $ \tau_j=v_{j}+\beta-\gamma-\frac{1}{p}+1, j=1, ..., n, \beta>\max\limits_{j}\left(\gamma+\frac{1}{p}-1-v_{j}\right).$

5)If $g\in M_T(A^{\infty, q}_\gamma, Y)$ then
\[\sup\limits_{r\in I^n}M_\infty(D^{\beta}g, r_1, ..., r_n)
\prod\limits_{j=1}^{n}(1-r_j)^{\tau_j}<\infty, \] $ \tau_j=v_{j}+\beta-\gamma-\frac{1}{p}+1, j=1, ..., n, \beta>\max\limits_{j}\left(\gamma+\frac{1}{p}-1-v_{j}\right).$
with $p=\infty$.
\end{theorem}
The following two theorems can be seen  as a corollary of lemma's 3 and 4 and also lemma 9 and 10 above. Lemma 3 and 4  provides estimates of Bergman kernel in  mentioned mixed norm   spaces $T^{p, q}_\gamma$,$M^{p, q}_\gamma$ including limit "infinity" cases of parameters and the closed graph theorem Note similar resuls are valid for $B^{p,q}_\alpha$ spaces.
We leave this to readers.To get estimates for Bergman kernel
for $B^{p,q}$ spaces we refer the reader to \cite{MitrShamoyan}.
We assume again $Y$ is a same BMOA space as in previous theorem.
\begin{theorem}
Let $p,q,\gamma>0$,$v>0$.
Let $g\in H(U^n); g(z)=\sum\limits_{k\in Z^n_{+}}c_kz^k$.  Then
\[1) If\  g\in (M_T)(M^{p, q}_\gamma, Y)\ then\ \sup\limits_{r\in I}M_\infty(D^\beta g, r)(1-r)^{v+n(\beta+1)-\gamma n-\frac{1}{p}}<\infty;\] for $\beta>\gamma+\frac{1}{np}-\frac{v}{n}-1$,
\[2) If\  g\in (M_T)(M^{p, \infty}_\gamma, Y)\ then\ \sup\limits_{r\in I}M_\infty(D^\beta g, r)(1-r)^{v+n(\beta+1)-\gamma n-\frac{1}{p}}<\infty;\]for $\beta>\gamma+\frac{1}{np}-\frac{v}{n}-1$,
\[3) If\  g\in (M_T)(M^{\infty, q}_\gamma, Y)\ then\ \sup\limits_{r\in I}M_\infty(D^\beta g, r)(1-r)^{v+n(\beta+1)-\gamma n}<\infty\]for $\beta>\gamma-\frac{v}{n}-1$,
\end{theorem}

Replacing $(M^{p, q}_\gamma); (M^{p, \infty}_\gamma); (M^{\infty, q}_\gamma)$ by $(T^{p, q}_\gamma); (T^{p, \infty}_\gamma); (T^{\infty, q}_\gamma)$ we have the following analogue of previous theorem for $T^{p,q}_{\alpha}$ type spaces of analytic function in polydisc based again on lemma's 3, 4 and 9,10 we provided above. Note these results though are not sharp but  are new even in case of unit disk.
We formulate this theorem the last one for $Y$ as 
$F^{\infty,d}_{\tau d-1,\alpha}$,for this case also $d$ parameter and $\alpha$ parameter do not appear in spaces of multipliers we found below.
We note that we consider general case when $\tau$ is a vector 
(see for this general case our discussion above)
\begin{theorem}
Let $p,q,\gamma>0$
Let  $g\in H(U^n), g(z)=\sum\limits_{k\in Z^n_{+}}c_kz^k$. Then
 \[1) If\  g\in (M_T)(T^{p, q}_\gamma, Y)\ then \ \sup\limits_{r\in I^n}M_\infty(D^\beta g, r)\prod\limits_{j=1}^{n}(1-r_j)^{\tau_j+\beta+1-\frac{\gamma}{n}-\frac{1}{p}}<\infty,\]\[ \beta>\max\limits_{j}\left(\frac{1}{p}+\frac{\gamma}{n}-\tau_j-1\right), \tau_j>0;\; j=1,\dots,n\]
 \[2) If\  g\in (M_T)(T^{p, \infty}_\gamma, Y)\ then\ \sup\limits_{r\in I^n}M_\infty(D^\beta g, r)\prod\limits_{j=1}^{n}(1-r_j)^{\tau_j+\beta+1-\frac{\gamma}{n}-\frac{1}{p}}<\infty, \]\[ \beta>\max\limits_{j}\left(\frac{1}{p}+\frac{\gamma}{n}-\tau_j-1\right), \tau_j >0;\; j=1,\dots,n\]
 \[3) If\  g\in (M_T)(T^{\infty, q}_\gamma, Y)\ then\ \sup\limits_{r\in I^n}M_\infty(D^\beta g, r)\prod\limits_{j=1}^{n}(1-r_j)^{\tau_j+\beta+1-\frac{\gamma}{n}-\frac{1}{p}}<\infty,\]\[ \beta>\max\limits_{j}\left(\frac{1}{p}+\frac{\gamma}{n}-\tau_j-1\right), \tau_j>0, \; j=1,\dots,n.\] 
with $p=\infty$.
\end{theorem}
We add another short comment here.
Note that we can also replace  $T^{p,q}_{\alpha}$ spaces in previous theorem by 
$H^{\vec{p}}_{\vec{\alpha}}$ classes .To formulate   such type new assertion  
we have to replace estimates of Bergman kernel of $T^{p,q}_{\alpha}$
spaces by estimates of Bergman kernel($g_{R}$ function) in quazinorms of $H^{\vec{p}}_{\vec{\alpha}}$ classes.Both can be seen in previous section. For $H^{\overrightarrow{p}}(\overrightarrow{\alpha})$  spaces  similar arguments can be easily applied.All results which are
needed for proofs  we just discussed for $H^{\overrightarrow{p}}(\overrightarrow{\alpha})$ and ,$H^{\vec{p}}_{\vec{\alpha}}$ analytic classes   can be found in previous section.We leave this procedure to for interested readers.

Short schemes of proofs of theorems 3, 4, 5.

Here is a short scheme of proof of various of type contained in theorem 3, 4, 5.

Assume $\{c_k\}_{k\in Z^n_+}\in M_T(X,Y)$ then we have to use the closed graph theorem and lemma 9,10 .

Indeed we have the following estimates and equalities, first by closed graph theorem $\|M_cf\|_{Y}\leq c\|f\|_X$. Let $w\in U^n$, and set $f_w(z)=\frac{1}{1-wz}$, $g_w=M_c(f_w)$. Then we have $$D^\beta g_w=D^\beta(M_c(f_w))=M_c(D^\beta(f_w))=cM_c\frac{1}{(1-wz)^{\beta+1}},\;\beta>0,$$ which gives $$\|D^\beta g_w\|_{Y}\leq c\left\|\frac{1}{(1-wz)^{\beta+1}}\right\|_X.$$
It remains to apply lemma 9,10 from left for BMOA type spaces from right for various $X$ spaces we have to use various estimates from lemma 3, 4, 6, 7 of Bergman kernel $g_R(z)$ in polydisc $U^n$ to get the desired result. 
We finnaly add one more remark.
\begin{remark}
Note probably some assertions of theorem 3 and 4 and 5 can be sharpened, namely the reverse implications can be proved also.
We do not discuss these issues in this paper.
\end{remark}
We can also use for some new assertions the following simple idea  (see \cite{Shvedenko}, \cite{JP}), \cite{MP1}) , the following obvious fact if $(g)\sim(c_k)$ is a multiplier from $X$ to $Y$ ($X, Y$ are quazynormed subspaces of $H(U^n)$) then it is a multiplier from dual spaces $Y^*$ to $X^*$ and then apply to this last pair a standard closed graph theorem. This  procedure was used before by many authors (see for example \cite{Shvedenko}).We can also use this type argument for our purposes.
Descriptions of duals of $H^{\overrightarrow{p}}(\overrightarrow{\alpha})$  spaces for $p_j>1, j=1, ..., n$ or $p_j<1, j=1, ..., n$ are known(see for this description \cite{Yaroslavcheva}, 
\cite{YaroslavchevaShamoyan}). 
Note also $(H^p)^*, (A^{\overline{p}, \overline{q}}_\alpha)^*, (A^{\widetilde{p}, \widetilde{q}}_\alpha)^*$ spaces are also well- known.(see ,for example, \cite{Shvedenko}) .
\footnotesize{}


\begin{thebibliography}{99}
\bibitem{Yaroslavcheva}
O.V. Yaroslavtseva, {\em Multipliers in some anisotropic spaces of analytic functions}, Investigations on linear operators and function theory. 26, Zap. Nauchn. Sem. POMI, 255, POMI, St. Petersburg., 1998, 244-248.
\bibitem{YaroslavchevaShamoyan}
F.A. Shamoyan, O.V. Yaroslavtseva, {\em Continuous projections, duality and the diagonal map in weighted spaces of holomorphic functions with mixed norm}, Investigations on linear operators and function theory. 25, Zap. Nauchn. Sem. POMI, 247, POMI, St. Petersburg., 1997, 268-275.
% \bibitem{ArcenovichShamoyan2007}
%M. Arsenovic, R. F. Shamoyan, {\em Sharp theorems on multipliers and distances in harmonic function spaces in higher dimension},  Zh. SFU. Ser. Mat. and phys., 5:3 (2012), 291-302.
\bibitem{Shvedenko}
S.V. Shvedenko, {\em Hardy classes and related spaces of analytic functions in the unit disc, polydisc and ball}, Results of Science and Technical. Ser. Mat. Anal., 23, VINITI, Moscow, 1985, 3-124.
\bibitem{ShamoyanMitrich}
 R. Shamoyan, O. Mihic, {\em On some inequalities in holomorphic function theory in polydisc related to diagonal mapping}, Chechoslovak Math. Journal, 2010.
\bibitem{Na1} M. Nawrocki, {\it Multipliers, linear functionals and the Frechet envelope of the Smirnov class
    $N_\ast(U^n)$}, Trans. Amer. Math. Soc., Vol 322, (1990), No. 2, 493-506.
\bibitem{AP} M.Arsenovic,R.Shamoyan , {\it On Fefferman-Stein type characterizations of certain analytic functional spaces
in unit disk},AMAPN,2012,v.4.
\bibitem{BP} O. Blasco, M. Pavlovic, {\it Coefficient multipliers in spaces of analytic functions}, preprint, 2011.
\bibitem{DS} M. Djrbashian and F. Shamoian, {\it Topics in the theory of $A^p_\alpha$ classes}, Teubner Texte zur
                       Mathematik, 1988, v 105.
\bibitem{GL} V. S. Guliev, P. I. Lizorkin, {\it Holomorphic and harmonic function classes in the polydisc and its boundary
    values}, in Russian, Trudi Math. Inst. RAN, Vol. 204, (1993), 137-159.
\bibitem{JP} M. Jevtic, M. Pavlovic, {\it Coefficient multipliers on spaces of analytic functions}, Acta Sci. Math., Vol. 64
   (1998), 531-545.
\bibitem{MP1} M. Matelievic, M. Pavlovic, {\it Multipliers of $H^p$ and $BMOA$}, Pacific J. Math., Vol. 146, (1990), 71-89.
\bibitem{OF1} J. Ortega, J. Fabrega, {\it Hardy's inequality and embeddings in holomorphic Lizorkin-Triebel spaces}, Illinois Journal of Mathematics, Vol. 43, 4, 1999, 733-751.
\bibitem{OF2} J. Ortega, J. Fabrega, {\it Pointwise multipliers and corona type decomposition in BMOA}, Ann. Inst. Fourier
    (Grenoble), Vol. 46, (1996), 111-137.
\bibitem{OF3} J. Ortega, J. Fabrega, {\it Holomoprhic Triebel-Lizorkin spaces}, J. Funct. Anal., Vol. 151, (1997), 177-212.
\bibitem{Sh1} R. F. Shamoyan, {\it On multipliers from Bergman type spaces into Hardy spaces in the polydisc}, in Russian,
   Ukrainian. Math. Journal., No. 10, (2000), 1405-1415.
\bibitem{Sh2} R. F. Shamoyan, {\it Continuous functionals and multipliers of a class of functions analytic in the polydisc},
   in Russian, Izv. Vuzov Matematika., No. 7, (2000), 67-69.
\bibitem{Sh3} R. F. Shamoyan, {\it Multipliers of power series Toeplitz operators and space embeddings}, Izv. Nat. Acad.
   Nauk Armenii. Vol. 34, (1999), No. 4, 43-59.
\bibitem{Tr1} R. M. Trigub, {\it Multipliers of the $H^p(U^m)$ class for $p \in (0,1)$, approximation properties of
    summability methods for power series}, in Russian, Doklady. RAN, Vol. 335, (1994), No. 6, 697-699.
\bibitem{ARS}R.Shamoyan,M.Arsenovic,{\it On multipliers of
mixed norm Lizorkin-Triebel type space},Journal of Siberian Federal University,2012,p.1-11.
\bibitem{MitrShamoyan}O. Mihic, R.Shamoyan, {\it On multipliers of some new analytic $M^{p, q}_\alpha$, $M^{p, \infty, \alpha}$ and $M^{\infty,p,\alpha}$ type spaces and related spaces on the unit polydisc}, Matematia Montisnigri, vol XXVII(2013) p.9-37.
\bibitem{KurShamoyan} R.Shamoyan,S.Kurilenko,{\it On multipliers
of some new Bergman-type and Hardy-type spaces in polydisk},
submitted,preprint,2013.
\bibitem{Shamoyan1} R.Shamoyan,{\it On multipliers of Bergman and Bloch spaces in the unit polydisk },Siberian Mathematical
Journal,2002.
\end{thebibliography}
\end{document}